\documentclass[a4paper,oneside]{amsart}

\RequirePackage{amsmath}
\RequirePackage{bm}
\RequirePackage{amssymb}
\RequirePackage{upref}
\RequirePackage{amsthm}
\RequirePackage{enumerate}
\RequirePackage{pb-diagram}
\RequirePackage{amsfonts}
\RequirePackage[mathscr]{eucal}
\RequirePackage{verbatim}
\RequirePackage{xr}
\RequirePackage{graphicx}
\RequirePackage[pdftex]{floatflt}
\usepackage{calc}
\usepackage{xspace}

\newcommand{\cf}{cf.\@\xspace}
\newcommand{\resp}{resp.\@\xspace}

\newcommand{\wlogc}{w.l.o.g.\@\xspace}
\newcommand{\Wlogc}{W.l.o.g.\@\xspace}

\newcommand{\al}{\alpha}
\newcommand{\bet}{\beta}
\newcommand{\ga}{\gamma}
\newcommand{\de}{\delta }
\newcommand{\e}{\epsilon}

\newcommand{\h}{\eta}

\newcommand{\ka}{\kappa}

\newcommand{\n}{\nu}
\newcommand{\om}{\omega}

\newcommand{\s}{\sigma}
\newcommand{\x}{\xi}

\newcommand{\C}{\varGamma}
\newcommand{\D}{\varDelta}

\newcommand{\Lam}{\varLambda}
\newcommand{\Om}{\varOmega}

\newcommand{\di}[1]{#1\nobreakdash-\hspace{0pt}dimensional}

\newcommand{\fv}[2]{#1\hspace{0pt}_{|_{#2}}}

\newcommand{\so}{{\mc S_0}}

\newcommand{\const}{\tup{const}}

\newcommand{\msp[1]}[1]{\mspace{#1mu}}

\newcommand{\R}[1][n+1]{{\protect\mathbb R}^{#1}}

\newcommand{\N}{{\protect\mathbb N}}

\newcommand{\eR}{\stackrel{\lower1ex \hbox{\rule{6.5pt}{0.5pt}}}{\msp[3]\R[]}}
\newcommand{\eN}{\stackrel{\lower1ex \hbox{\rule{6.5pt}{0.5pt}}}{\msp[1]\N}}
\newcommand{\eO}{\stackrel{\lower1ex
\hbox{\rule{6pt}{0.5pt}}}{\msc O}}

\DeclareMathOperator{\graph}{graph}

\newcommand\ra{\rightarrow}

\newcommand\pa{\partial}
\newcommand\pde[2]{\frac {\partial#1}{\partial#2}}
 
\newcommand{\un}{\infty}
\newcommand{\A}{\forall}

\newcommand{\set}[2]{\{\,#1\colon #2\,\}}
\newcommand{\uu}{\cup}

\newcommand{\uuu}{\bigcup}

\newcommand{\uud}{ \stackrel{\lower 1ex \hbox {.}}{\uu}}
\newcommand{\uuud}[1]{ \stackrel{\lower 1ex \hbox {.}}{\uuu_{#1}}}
\newcommand\su{\subset}

\newcommand{\sminus}[1][28]{\raise 0.#1ex\hbox{$\scriptstyle\setminus$}}

\newcommand{\wed}{\wedge}

\newcommand{\abs}[1]{\lvert#1\rvert}

\newcommand{\norm}[1]{\lVert#1\rVert}

\newcommand{\spd}[2]{\protect\langle #1,#2\protect\rangle}

\newcommand\ch[3]{\varGamma_{#1#2}^#3}
\newcommand\cha[3]{{\bar\varGamma}_{#1#2}^#3}

\newcommand{\riem}[4]{R_{#1#2#3#4}}
\newcommand{\riema}[4]{{\bar R}_{#1#2#3#4}}

\newcommand{\tit}{\textit}

\newcommand{\tup}{\textup}

\newcommand{\mc}{\protect\mathcal}
\newcommand{\msc}{\protect\mathscr}

\providecommand{\bysame}{\makebox[3em]{\hrulefill}\thinspace}

\newcommand{\ci}{\cite}

\newcommand{\bt}{\begin{thm}}
\newcommand{\bl}{\begin{lem}}
\newcommand{\bc}{\begin{cor}}
\newcommand{\bd}{\begin{definition}}
\newcommand{\bpp}{\begin{prop}}
\newcommand{\br}{\begin{rem}}
\newcommand{\bn}{\begin{note}}
\newcommand{\be}{\begin{ex}}
\newcommand{\bes}{\begin{exs}}
\newcommand{\bb}{\begin{example}}
\newcommand{\bbs}{\begin{examples}}
\newcommand{\ba}{\begin{axiom}}

\newcommand{\et}{\end{thm}}
\newcommand{\el}{\end{lem}}
\newcommand{\ec}{\end{cor}}
\newcommand{\ed}{\end{definition}}
\newcommand{\epp}{\end{prop}}
\newcommand{\er}{\end{rem}}
\newcommand{\en}{\end{note}}
\newcommand{\ee}{\end{ex}}
\newcommand{\ees}{\end{exs}}
\newcommand{\eb}{\end{example}}
\newcommand{\ebs}{\end{examples}}
\newcommand{\ea}{\end{axiom}}

\newcommand{\bp}{\begin{proof}}
\newcommand{\ep}{\end{proof}}
\newcommand{\eps}{\renewcommand{\qed}{}\end{proof}}

\newcommand{\bal}{\begin{align}}

\newcommand{\bi}[1][1.]{\begin{enumerate}[\upshape #1]}
\newcommand{\bia}[1][(1)]{\begin{enumerate}[\upshape #1]}
\newcommand{\bin}[1][1]{\begin{enumerate}[\upshape\bfseries #1]}
\newcommand{\bir}[1][(i)]{\begin{enumerate}[\upshape #1]}
\newcommand{\bic}[1][(i)]{\begin{enumerate}[\upshape\hspace{2\cma}#1]}
\newcommand{\bis}[2][1.]{\begin{enumerate}[\upshape\hspace{#2\parindent}#1]}
\newcommand{\ei}{\end{enumerate}}

\newcommand\ndots{\raise 0.47ex \hbox {,}\hskip0.06em\cdots %
     \raise 0.47ex \hbox {,}\hskip0.06em} 

\newcommand{\q}{\quad}
\newcommand{\qq}{\qquad}

\newcommand{\hp}{\hphantom}

\newcommand\nd{\noindent}

\newskip\Csmallskipamount                                                
\Csmallskipamount=\smallskipamount
\newskip\Cmedskipamount
\Cmedskipamount=\medskipamount
\newskip\Cbigskipamount
\Cbigskipamount=\bigskipamount

\newcommand\cvs{\vspace\Csmallskipamount}   
\newcommand\cvm{\vspace\Cmedskipamount}

\newskip\csa
\csa=\smallskipamount

\newskip\cma
\cma=\medskipamount

\newskip\cba
\cba=\bigskipamount

\newdimen\spt
\newcommand\citem{\cvs\advance\itemno by
1{(\romannumeral\the\itemno})\hskip3pt}
\newcommand{\bitem}{\cvm\nd\advance\itemno by
1{\bf\the\itemno}\hspace{\cma}}

\newcount\itemno
\newcommand{\lae}[1]{\label{E:#1}}
\newcommand{\lat}[1]{\label{T:#1}}

\newcommand{\lad}[1]{\label{D:#1}}

\newcommand{\rt}[1]{Theorem~\ref{T:#1}}

\newcommand{\rd}[1]{Definition~\ref{D:#1}}

\newcommand{\re}[1]{\eqref{E:#1}}

\newskip\thmskip
\thmskip=\parindent

\newskip\hsk
\newenvironment{hinw}{\labelsep=0pt\begin{list}{}{\labelsep=0pt\itemindent=0pt\labelwidth=0pt\leftmargin=\parindent\rightmargin=0pt\partopsep=\cba}%
\item\it\nopagebreak\nopagebreak}%
{\end{list}}

\newcommand\bh{\begin{hinw}}
\newcommand{\eh}{\end{hinw}}

\newtheoremstyle{normal}
  {\cba}
  {\cba}
  {}
  {\thmskip}
  {\bfseries}
  {.}
  {\hsk}
  {}

\newtheoremstyle{abschnitt}
  {\cba}
  {\cba}
  {}
  {\thmskip}
  {\bfseries}
  {.}
  {\hsk}
  {}

\newtheoremstyle{italic}
  {\cba}
  {\cba}
  {\itshape}
  {\thmskip}
  {\bfseries}
  {.}
  {\hsk}
  {}

\newtheoremstyle{aufgaben}
  {\cba}
  {\cba}
  {}
  {}
  {\normalsize\bfseries}
  {.}
  {\hsk}
  {}

\newtheoremstyle{break}
  {\cba}
  {\cba}
  {\itshape}
  {}
  {\bfseries}
  {.}
  {\newline}
  {}

\swapnumbers
\theoremstyle{italic}
\newtheorem{thm}[subsection]{Theorem}
\newtheorem{lem}[subsection]{Lemma}
\newtheorem{prop}[subsection]{Proposition}
\newtheorem{cor}[subsection]{Corollary}

\theoremstyle{normal}
\newtheorem{rem}[subsection]{Remark}
\newtheorem{definition}[subsection]{Definition}
\newtheorem{example}[subsection]{Example}
\newtheorem{examples}[subsection]{Examples}
\newtheorem{ex}[subsection]{Exercise}
\newtheorem{note}[subsection]{}
\newtheorem{axiom}[subsection]{Axiom}

\theoremstyle{aufgaben}
\newtheorem{exs}[subsection]{Exercises}

\swapnumbers

\numberwithin{equation}{section}
\numberwithin{figure}{section}

\newenvironment{textequation}[1][0.8]
{\begin{equation}
\begin{aligned}
\begin{minipage}{#1\linewidth}}
{\end{minipage}
\end{aligned}
\end{equation}
\ignorespacesafterend}

\newcommand{\btext}{\begin{textequation}}
\newcommand{\etext}{\end{textequation}}

\usepackage[german,english]{babel}
\usepackage{graphicx}

\RequirePackage{upref}
\RequirePackage{amsthm}
\RequirePackage{enumerate}
\usepackage[mathscr]{eucal}

\usepackage{xr-hyper}
\usepackage{hyperref}

\usepackage{calc}

\newlength{\oddsidemarginlength}
\newlength{\topmarginlength}

\newcounter{numberoflines}
\newcounter{tempcc}
\oddsidemargin=\oddsidemarginlength
\evensidemargin=\oddsidemargin
\topmargin=\topmarginlength
\begin{document}

\flushbottom


\title{The mass of a Lorentzian manifold}

\author{Claus Gerhardt}
\address{Ruprecht-Karls-Universit\"at, Institut f\"ur Angewandte Mathematik,
Im Neuenheimer Feld 294, 69120 Heidelberg, Germany}
\email{gerhardt@math.uni-heidelberg.de}
\urladdr{http://www.math.uni-heidelberg.de/studinfo/gerhardt/}
\thanks{This work has been supported by the Deutsche Forschungsgemeinschaft}

%
\date{\today}
%


\begin{abstract}
We define a physically reasonable mass for an asymptotically Robertson-Walker (ARW)
mani\-fold which is uniquely defined in the case of a normalized representation.
\end{abstract}

\maketitle

\tableofcontents
\setcounter{section}{0}
\section{Introduction}

For asymptotically flat Lorentzian 4-manifolds the so-called ADM-mass is defined by
looking at a spacelike slice $M=\{t=\const\}$. If $M$ is asymptotically flat and the
scalar curvature of $M$ of class $L^1(M)$, then the ADM-mass is defined as a flux
integral

\begin{equation}
m_{ADM}=\lim_{r\ra\un}\frac1{16\pi}\int_{\pa B_r(0)}(g_{ii,j}-g_{ij,i})\nu^j;
\end{equation}
the limit is finite iff $\int_M\abs R<\un$, \cf  \cite{br:mass}. Schoen and Yau
\cite{sy:mass:one, sy:mass} proved that $m\ge 0$ if $R\ge 0$ with equality if and only
if
$M$ is isometric to Euclidean space.

If $M$ is an exterior region with compact boundary consisting of minimal surfaces,
then the Penrose inequality states that

\begin{equation}
16\pi m^2\ge \abs C,
\end{equation}
where $\abs C$ is the area of any connected component of $\pa M$. The Penrose
inequality was proved by Huisken and Ilmanen \cite{hi:penrose}.

\cvm
For general Lorentzian manifolds, e.g., for cosmological spactimes there is no notion
of mass as far as we know. We shall define a physically reasonable mass for
asymptotically Robertson-Walker (ARW) spacetimes---see \rd{1.2}---that satisfy the
timelike convergence condition
\begin{equation}
\bar R_{\al\bet}\nu^\al\nu^\bet \ge 0\qq\A\,\spd\n\nu=-1.
\end{equation}

\cvm
In \cite{cg:arw} we introduced the notion of asymptotically Robertson-Walker
spacetimes and proved  some convergence theorems for solutions to the inverse
mean curvature flow (IMCF) in these spacetimes.

\bd\lad{1.2}
A cosmological spacetime $N$, $\dim N=n+1$, is said to be \tit{asymptotically
Robertson-Walker} (ARW) with respect to the future, if a future end of $N$, $N_+$,
can be written as a product $N_+=[a,b)\times \so$, where $\so$ is a compact
Riemannian space, and there exists a future directed time function $\tau=x^0$ such
that the metric in $N_+$ can be written as
\begin{equation}\lae{1.19}
d\bar s^2=e^{2\tilde\psi}\{-{dx^0}^2+\s_{ij}(x^0,x)dx^idx^j\},
\end{equation}
where  $\so$ corresponds to $x^0=a$, $\tilde\psi$ is of the form
\begin{equation}
\tilde\psi(x^0,x)=f(x^0)+\psi(x^0,x),
\end{equation}
and we assume that there exists a positive constant $c_0$ and a smooth
Riemannian metric $\bar\s_{ij}$ on $\so$ such that
\begin{equation}
\lim_{\tau\ra b}e^\psi=c_0\q\wed\q \lim_{\tau\ra b}\s_{ij}(\tau,x)=\bar\s_{ij}(x),
\end{equation}
and
\begin{equation}
\lim_{\tau\ra b}f(\tau)=-\un.
\end{equation}

\cvm
Without loss of generality we shall assume $c_0=1$. Then $N$ is ARW with
respect to the future, if the metric is close to the Robertson-Walker metric
\begin{equation}\lae{0.8}
d\bar s^2=e^{2f}\{-{dx^0}^2+\bar\s_{ij}(x)dx^idx^j\}
\end{equation}
near the singularity $\tau =b$. By \tit{close} we mean that the derivatives of arbitrary order with respect to space and time of the
conformal metric $e^{-2f}\bar g_{\al\bet}$ in \re{1.19} should converge  to the corresponding derivatives of the conformal limit metric in \re{0.8} when $x^0$ tends to $b$. We
emphasize that in our terminology Robertson-Walker metric does not imply that
$(\bar\s_{ij})$ is a metric of constant curvature, it is only the spatial metric of a
warped product.

\cvm
We assume, furthermore, that $f$ satisfies the following five conditions
\begin{equation}
-f'>0,
\end{equation}
there exists $\om\in\R[]$ such that
\begin{equation}\lae{1.25}
n+\om-2>0\q\wed\q \lim_{\tau\ra b}\abs{f'}^2e^{(n+\om-2)f}=m>0.
\end{equation}
Set $\tilde\ga =\frac12(n+\om-2)$, then there exists the limit
\begin{equation}\lae{0.11}
\lim_{\tau\ra b}(f''+\tilde\ga \abs{f'}^2)
\end{equation}
and
\begin{equation}
\abs{D^m_\tau(f''+\tilde\ga \abs{f'}^2)}\le c_m \abs{f'}^m\qq
\A\, m\ge 1,
\end{equation}
as well as
\begin{equation}
\abs{D_\tau^mf}\le c_m \abs{f'}^m\qq\A\, m\ge 1.
\end{equation}

\cvm
We call $N$ a \tit{normalized} ARW spacetime, if
\begin{equation}
\int_{\so}\sqrt{\det{\bar\s_{ij}}}=\abs{S^n}.
\end{equation}
\ed

\br
(i) If these assumptions are satisfied then the range of $\tau$ is finite, hence, we
may---and shall---assume \wlogc that $b=0$, i.e.,
\begin{equation}
a<\tau<0.
\end{equation}

\cvm
(ii) Any ARW spacetime can be normalized as one easily checks. Without a
normalization condition the constant $m$ in \re{1.25} wouldn't be defined uniquely
as we shall see. It will later be identified with the mass of $N$.

\cvm
(iii) In view of the assumptions on $f$ the mean curvature of the coordinate slices
$M_\tau=\{x^0=\tau\}$ tends to $\un$, if $\tau$ goes to zero.

\cvm
(iv) Similarly one can define $N$ to be ARW with respect to the past. In this case the
singularity would lie in the past, correspond to $\tau=0$, and the mean curvature
of the coordinate slices would tend to $-\un$.
\er

Our main result is

\bt\lat{0.3}
Let $N$ be a $(n+1)$-dimensional normalized  ARW spacetime with respect to the
future  that satisfies the timelike convergence condition.
Then the future mass $m$ of $N$ is defined by
\begin{equation}
\tfrac12 n(n-1)\abs{S^n}m=\lim \int_M G_{\al\bet}\nu^\al\nu^\bet e^{\om f}
e^{\psi},
\end{equation}
where $G_{\al\bet}$ is the Einstein tensor, $\om$ the constant that appears in the
definition of ARW spaces, and the closed spacelike hypersurfaces
$M$ converge to the future singularity such that, if they are written as graphs over
$\so$, $M=\graph u$, $Du$ vanishes when the singularity is approached.

For normalized ARW spaces the mass is defined independently of the time function
$x^0$, $f$, $\psi$ and $\om$.
\et

\section{Notations and definitions}

The main objective of this section is to state the equations of Gau{\ss}, Codazzi,
and Weingarten for space-like hypersurfaces $M$ in a \di {(n+1)} Lorentzian
manifold
$N$.  Geometric quantities in $N$ will be denoted by
$(\bar g_{ \al \bet}),(\riema  \al \bet \ga \de)$, etc., and those in $M$ by $(g_{ij}), 
(\riem ijkl)$, etc.. Greek indices range from $0$ to $n$ and Latin from $1$ to $n$;
the summation convention is always used. Generic coordinate systems in $N$ resp.
$M$ will be denoted by $(x^ \al)$ resp. $(\x^i)$. Covariant differentiation will
simply be indicated by indices, only in case of possible ambiguity they will be
preceded by a semicolon, i.e., for a function $u$ in $N$, $(u_ \al)$ will be the
gradient and
$(u_{ \al \bet})$ the Hessian, but e.g., the covariant derivative of the curvature
tensor will be abbreviated by $\riema  \al \bet \ga{ \de;\e}$. We also point out that
\begin{equation}
\riema  \al \bet \ga{ \de;i}=\riema  \al \bet \ga{ \de;\e}x_i^\e
\end{equation}
with obvious generalizations to other quantities.

Let $M$ be a \tit{spacelike} hypersurface, i.e., the induced metric is Riemannian,
with a differentiable normal $\n$ which is time-like.

In local coordinates, $(x^ \al)$ and $(\x^i)$, the geometric quantities of the
space-like hypersurface $M$ are connected through the following equations
\begin{equation}\lae{1.2}
x_{ij}^ \al= h_{ij}\n^ \al
\end{equation}
the so-called \tit{Gau{\ss} formula}. Here, and also in the sequel, a covariant
derivative is always a \tit{full} tensor, i.e.

\begin{equation}
x_{ij}^ \al=x_{,ij}^ \al-\ch ijk x_k^ \al+ \cha  \bet \ga \al x_i^ \bet x_j^ \ga.
\end{equation}
The comma indicates ordinary partial derivatives.

In this implicit definition the \tit{second fundamental form} $(h_{ij})$ is taken
with respect to $\n$.

The second equation is the \tit{Weingarten equation}
\begin{equation}
\n_i^ \al=h_i^k x_k^ \al,
\end{equation}
where we remember that $\n_i^ \al$ is a full tensor.

Finally, we have the \tit{Codazzi equation}
\begin{equation}
h_{ij;k}-h_{ik;j}=\riema \al \bet \ga \de\n^ \al x_i^ \bet x_j^ \ga x_k^ \de
\end{equation}
and the \tit{Gau{\ss} equation}
\begin{equation}
\riem ijkl=- \{h_{ik}h_{jl}-h_{il}h_{jk}\} + \riema  \al \bet\ga \de x_i^ \al x_j^ \bet
x_k^ \ga x_l^ \de.
\end{equation}

Now, let us assume that $N$ is a globally hyperbolic Lorentzian manifold with a
\tit{compact} Cauchy surface. 
$N$ is then a topological product $I\times \mc S_0$, where $I$ is an open interval,
$\mc S_0$ is a compact Riemannian manifold, and there exists a Gaussian coordinate
system
$(x^ \al)$, such that the metric in $N$ has the form 
\begin{equation}\lae{1.7}
d\bar s_N^2=e^{2\psi}\{-{dx^0}^2+\s_{ij}(x^0,x)dx^idx^j\},
\end{equation}
where $\s_{ij}$ is a Riemannian metric, $\psi$ a function on $N$, and $x$ an
abbreviation for the spacelike components $(x^i)$. 
We also assume that
the coordinate system is \tit{future oriented}, i.e., the time coordinate $x^0$
increases on future directed curves. Hence, the \tit{contravariant} time-like
vector $(\x^ \al)=(1,0,\dotsc,0)$ is future directed as is its \tit{covariant} version
$(\x_ \al)=e^{2\psi}(-1,0,\dotsc,0)$.

Let $M=\graph \fv u\so$ be a space-like hypersurface
\begin{equation}
M=\set{(x^0,x)}{x^0=u(x),\,x\in\mc S_0},
\end{equation}
then the induced metric has the form
\begin{equation}
g_{ij}=e^{2\psi}\{-u_iu_j+\s_{ij}\}
\end{equation}
where $\s_{ij}$ is evaluated at $(u,x)$, and its inverse $(g^{ij})=(g_{ij})^{-1}$ can
be expressed as
\begin{equation}\lae{1.10}
g^{ij}=e^{-2\psi}\{\s^{ij}+\frac{u^i}{v}\frac{u^j}{v}\},
\end{equation}
where $(\s^{ij})=(\s_{ij})^{-1}$ and
\begin{equation}\lae{1.11}
\begin{aligned}
u^i&=\s^{ij}u_j\\
v^2&=1-\s^{ij}u_iu_j\equiv 1-\abs{Du}^2.
\end{aligned}
\end{equation}
Hence, $\graph u$ is space-like if and only if $\abs{Du}<1$.

The covariant form of a normal vector of a graph looks like
\begin{equation}
(\n_ \al)=\pm v^{-1}e^{\psi}(1, -u_i).
\end{equation}
and the contravariant version is
\begin{equation}
(\n^ \al)=\mp v^{-1}e^{-\psi}(1, u^i).
\end{equation}
Thus, we have
\br Let $M$ be space-like graph in a future oriented coordinate system. Then the
contravariant future directed normal vector has the form
\begin{equation}
(\n^ \al)=v^{-1}e^{-\psi}(1, u^i)
\end{equation}
and the past directed
\begin{equation}\lae{1.15}
(\n^ \al)=-v^{-1}e^{-\psi}(1, u^i).
\end{equation}
\er

In the Gau{\ss} formula \re{1.2} we are free to choose the future or past directed
normal, but we stipulate that we always use the past directed normal for reasons
that we have explained in \ci[Section 2]{cg:indiana}.

Look at the component $ \al=0$ in \re{1.2} and obtain in view of \re{1.15}

\begin{equation}\lae{1.16}
e^{-\psi}v^{-1}h_{ij}=-u_{ij}- \cha 000\mspace{1mu}u_iu_j- \cha 0j0
\mspace{1mu}u_i- \cha 0i0\mspace{1mu}u_j- \cha ij0.
\end{equation}
Here, the covariant derivatives are taken with respect to the induced metric of
$M$, and
\begin{equation}
-\cha ij0=e^{-\psi}\bar h_{ij},
\end{equation}
where $(\bar h_{ij})$ is the second fundamental form of the hypersurfaces
$\{x^0=\const\}$.

An easy calculation shows
\begin{equation}
\bar h_{ij}e^{-\psi}=-\tfrac{1}{2}\dot\s_{ij} -\dot\psi\s_{ij},
\end{equation}
where the dot indicates differentiation with respect to $x^0$.

\section{Proof of \rt{0.3}}

Let $N$ be a normalized ARW spacetime and assume that the metric is given by
\re{1.19} such that the future end is described by
\begin{equation}
a\le x^0<0.
\end{equation}
\Wlogc we also suppose
\begin{equation}
\lim_{x^0\ra 0}\psi(x^0,x)=0.
\end{equation}

Let $(\tilde g_{\al\bet})$ be the conformal metric
\begin{equation}\lae{2.3}
\tilde g_{\al\bet}dx^\al dx^\bet= -(dx^0)^2+\s_{ij}(x^0,x)dx^idx^j
\end{equation}
and distinguish the geometric quantities with respect to this metric by a tilde, i.e.,
$\tilde R_{\al\bet\ga\de}$ is the Riemannian curvature tensor, etc.. Then we have
\begin{equation}\lae{2.4}
\bar R_{\al\bet}=\tilde
R_{\al\bet}-(n-1)[\tilde\psi_{\al\bet}-\tilde\psi_\al\tilde\psi_\bet]-\tilde
g_{\al\bet} [\D\tilde\psi + (n-1) \norm{D\tilde\psi}^2],
\end{equation}
and
\begin{equation}
\bar R=e^{-2\tilde\psi}[\tilde R -2n\D\tilde\psi -n(n-1)\norm{D\tilde\psi}^2],
\end{equation}
where the covariant derivatives of $\tilde\psi$ are taken with respect to the metric
$(\tilde g_{\al\bet})$.

\cvm
The Einstein tensor in $N$ is defined by
\begin{equation}
G_{\al\bet}=\bar R_{\al\bet}-\tfrac12 \bar R \bar g_{\al\bet}
\end{equation}
It is divergent free, i.e.,
\begin{equation}
G^\al_{\bet;\al}=0.
\end{equation}

\cvm
Let $(\h_\al)=e^{\tilde\psi}(-1,0,\ldots,0)$ be the covariant vector field that
represents the future directed normal of the slices $\{x^0=\const\}$ and let
$\Om\su N$ be  an open subset bounded by two spacelike hypersurfaces $M_1$ and
$M_2$, where $M_2$ should lie in the future of $M_1$. Applying Gau{\ss}'
divergence theorem we obtain
\begin{equation}\lae{2.8}
\begin{aligned}
0=&\int_\Om G^{\al\bet}_{\hp{\al\bet};\al}\h_\bet e^{\om f}e^\psi\\[\cma]
=&-\int_\Om G^{\al\bet}\h_{\bet;\al}e^{\om f} e^\psi-\int_\Om G^{\al\bet}
\h_\bet [\om f_\al+\psi_\al]e^{\om f} e^\psi\\[\cma]
&-\big(\int_{M_1}G_{\al\bet}\nu^\al \h^\bet e^{\om f} e^\psi
+\int_{M_2}G_{\al\bet}
\nu^\al \h^\bet e^{\om f} e^\psi\big),
\end{aligned}
\end{equation}
where the normals $\nu$ of the hypersurfaces $M_i$ are supposed to point
outward of $\Om$, i.e., in case of $M_1$, $\nu$ is past directed, and in case of
$M_2$, $\nu$ is future directed. Note the minus sign in front of the boundary
integrals which is due to the sign of
$\spd\nu\nu=-1$.

The covariant derivatives of $(\h_\al)$, $(\h_{\al;\bet})$, satisfy
\begin{equation}
\h_{0;\bet}=0,\qq\h_{i;0}=\bar\C^0_{i0} e^{\tilde\psi}=\psi_{,i}e^{\tilde\psi},
\end{equation}
where $\psi_{,i}=\pde \psi {x^i}$, and
\begin{equation}
\h_{i;j}=\bar\C^0_{ij} e^{\tilde\psi}=-\bar h_{ij},
\end{equation}
where $(\bar h_{ij})$ is the second fundamental form of the slices $\{x^0=\const\}$.

Hence we deduce from \re{2.8}
\begin{equation}\lae{2.11}
\begin{aligned}
&\int_{M_1}G_{\al\bet}\nu^\al \h^\bet e^{\om f} e^\psi
+\int_{M_2}G_{\al\bet}
\nu^\al \h^\bet e^{\om f} e^\psi\\[\cma]
=&\int_\Om G^{ij}\bar h_{ij} e^{\om f} e^\psi + \int_\Om G^{00}[\om f'+\psi']
e^{\om f} e^{\tilde\psi}e^\psi\\[\cma]
=&\int_\Om n(n-1)[f''+\tilde\ga \abs{f'}^2]f' e^{(\om-3)f} \\[\cma]
&+\int_\Om c\, e^{(\om-3)f},
\end{aligned}
\end{equation}
where the symbol $c$ represents  terms that can be estimated by
\begin{equation}
\abs c\le c_0(1+\abs{f'}+\e\abs{f'}^2)\q\tup{and}\q\lim_{x^0\ra 0}\e=0.
\end{equation}
To derive the second equality in \re{2.11} we used the relations \re{2.3}, \re{2.4}
as well as the assumption that the metrics $\s_{ij}(\tau,\cdot)$ converge in
$C^\un$ to $\bar\s_{ij}$.

The volume element in $\Om$ is of the form
\begin{equation}
e^{(n+1)f} e^{(n+1)\psi} \sqrt{\det(\s_{ij})}\,dx dx^0,
\end{equation}
thus the right-hand side of \re{2.11} vanishes if the hypersurfaces $M_i$ approach
the singularity, in view of \re{1.25} and \re{0.11}.

\cvm
Now, in \re{2.11} let us choose the hypersurfaces $M_i$ to be  slices
$\{x^0=\const\}$, then the left-hand side is equal to
\begin{equation}
\int_{M_2}G_{\al\bet}\nu^\al \nu^\bet e^{\om f} e^\psi
-\int_{M_1}G_{\al\bet}
\nu^\al \nu^\bet e^{\om f} e^\psi
\end{equation}
and we conclude that
\begin{equation}
\lim_{\tau\ra 0}\int_{M_\tau}G_{\al\bet}\nu^\al \nu^\bet e^{\om f} e^\psi
\end{equation}
exists, where $M_\tau=\{x^0=\tau\}$, and the limit is equal to
\begin{equation}
\tfrac12 n(n-1)\lim_{\tau\ra 0}\abs{f'}^2 e^{(n+\om -2)f}\int_\so
\sqrt{\det(\bar\s_{ij})}=\tfrac12 n(n-1) m\abs{S^n}.
\end{equation}

From the above considerations we immediately deduce that the preceding limit is
equal to
\begin{equation}
\lim \int_{M_k} G_{\al\bet}\nu^\al\nu^\bet e^{\om f} e^\psi,
\end{equation}
where $M_k=\graph u_k$ are arbitrary spacelike hypersurfaces, written as graphs
over $\so$, such that
\begin{equation}
\lim u_k=0\q \tup{and}\q \lim \abs{Du_k}=0.
\end{equation}

\cvm
Hence we may use the leaves $M(t)$ of an inverse mean curvature flow
\begin{equation}\lae{2.18}
\dot x=-H^{-1}\nu
\end{equation}
with initial hypersurface $M_0$, $\fv H{M_0}>0$, to define the mass, since the flow
hypersurfaces $M(t)$ run straight into the singularity and satisfy
\begin{equation}\lae{2.19}
\abs u_m\le c_m e^{-\ga t}\qq \A\, m\in \N,
\end{equation}
where $\ga=\frac1n \tilde\ga$, \cf \cite[Lemma 7.1]{cg:arw}.

\cvm
Using the Gau{\ss} equation
\begin{equation}
R=-[H^2-\abs{A}^2]+2G_{\al\bet}\nu^\al\nu^\bet
\end{equation}
we can rewrite
\begin{equation}\lae{2.21}
\begin{aligned}
\int_{M} G_{\al\bet}\nu^\al\nu^\bet e^{\om f} e^\psi &=\tfrac{n-1}{2n} \int_M
H^2 e^{\om f}e^\psi \\[\cma]
&\hp{=}+ \tfrac12 \int_M\big(R-[\abs A^2-\tfrac1n H^2]\big) e^{\om f}
e^\psi
\end{aligned}
\end{equation}
to conclude that
\begin{equation}
\lim \int_{M} G_{\al\bet}\nu^\al\nu^\bet e^{\om f} e^\psi=\lim \tfrac{n-1}{2n} \int_M
H^2 e^{\om f}e^\psi 
\end{equation}
for those hypersurfaces for which the second integral on the right-hand side of
\re{2.21} tends to zero if  the singularity is approached.

\cvm
This is the case for the coordinate slices $\{x^0=\const\}$ as well as for the leaves
of an IMCF.

\bl
Let $M(t)$ be a solution of the evolution equation \re{2.18}, then
\begin{equation}
\lim_{t\ra \un}\int_{M(t)}\big(R-[\abs A^2-\tfrac1n H^2]\big)e^{\om f}e^\psi =0.
\end{equation}
\el

\bp
(i) Let us first estimate the scalar curvature. We have
\begin{equation}
\bar g_{\al\bet}=e^{2\tilde\psi}\tilde g_{\al\bet},
\end{equation}
where $(\tilde g_{\al\bet})$ is the metric in \re{2.3}. Denote by $\tilde h_{ij}, \tilde
g_{ij}$, etc. the geometric quantities of hypersurfaces when the metric of the
ambient space is $(\tilde g_{\al\bet})$, then
\begin{align}
e^{\tilde\psi}h^j_i&=\tilde h^j_i+\tilde\psi_\al \tilde\nu^\al
\de^j_i,\lae{2.25}\\[\cma] 
g_{ij}&=e^{2\tilde\psi}\tilde
g_{ij}=e^{2\tilde\psi}(-u_iu_j+\s_{ij}(u,x)dx^idx^j),\\
\intertext{and}
R&=e^{-2\tilde\psi}(\tilde R-2(n-1)\D \tilde\psi-(n-1)(n-2)\norm{D\tilde\psi}^2),
\end{align}
where the covariant derivatives of $\tilde\psi(u,x)$ are taken with respect to $\tilde
g_{ij}$.

\cvm
Now, $\tilde R$ is bounded and the covariant derivatives of $\tilde \psi$ are bounded
as well, \cf \cite[Section 6]{cg:arw}, hence
\begin{equation}
\int_M\abs R e^{\om f}e^\psi\le c\int_\so e^{(n+\om-2)f}
e^{(n+1)\psi}\sqrt{\det(\s_{ij})}\ra 0.
\end{equation}

\cvm
(ii) To estimate the second fundamental form, we use \re{2.25} to obtain
\begin{equation}
e^{2\tilde \psi}(\abs A^2-\tfrac1n H^2)=\abs{\tilde A}^2-\tfrac1n \tilde H^2
\end{equation}
and thus
\begin{equation}
\begin{aligned}
&\int_M(\abs A^2-\tfrac1n H^2)e^{\om f} e^\psi\\[\cma]
&=\int_\so (\abs{\tilde
A}^2-\tfrac1n \tilde H^2)e^{(n+\om-2)f}e^{(n-1)\psi}\sqrt{\det(\s_{ij})}\ra 0.
\end{aligned}
\end{equation}
\ep

The final part of \rt{0.3} is proved in the next section.

\section{Uniqueness of the mass}

The mass of a normalized ARW space doesn't depend on the particular time function.
For a proof of this claim we shall once again employ the leaves of an IMCF.

\cvm
Let $\tilde x^0$ be a second time function that provides a normalized representation
of $N$ as an ARW space
\begin{equation}
d\bar s^2=e^{2(\tilde f+\tilde\psi)}(-(d\tilde x^0)^2+\tilde \s_{ij}dx^idx^j)
\end{equation}
such that
\begin{equation}\lae{3.2}
\lim\abs{\tilde f'}^2e^{(n+\tilde\om -2)\tilde f}=\tilde m>0.
\end{equation}

\cvm
Let $M(t)$ be a solution of the evolution problem \re{2.18} which are written as
$M(t)=\graph u$ in the original coordinate system and as $M(t)=\graph \tilde u$ in
the new system. In both cases we may assume that the hypersurfaces are graphs
over $\so$, since it is unnecessary that $\so$ is a level hypersurface for a time
function. The estimates \re{2.19} are satisfied by $u$ as well as $\tilde u$.

There are many invariants that could be used to compare $f$ and $\tilde f$. Let us
consider $G_{\al\bet}\nu^\al\nu^\bet$ evaluated at $M(t)$. Arguing as in the
preceding section we deduce
\begin{equation}
1=\lim_{t\ra\un}\frac{\abs{f'}^2e^{(n+\om-2)f}}{\abs{\tilde f'}^2e^{(n+\om-2)
\tilde f}} \frac{e^{(n+\tilde\om)\tilde f}}{e^{(n+\om)f}},
\end{equation}
where the arguments of $f, f'$ \resp $\tilde f, \tilde f'$ are $u$ \resp $\tilde u$.

Hence we conclude
\begin{equation}\lae{3.4}
\lim \frac{e^{(n+\tilde\om)\tilde f}}{e^{(n+\om)f}}=\frac{\tilde m}m,
\end{equation}
in view of \re{1.25} and \re{3.2}.

\cvm
Now we observe that
\begin{equation}
\frac d{dt} f(u)=f'\dot u=\frac{\tilde v f'}{\tilde H-n\tilde v f'+\psi_\al\tilde\nu^\al},
\end{equation}
where $\tilde H$ is the mean curvature with respect to the conformal metric in
\re{2.3} and $\tilde v=v^{-1}$. Hence we obtain
\begin{equation}\lae{3.6}
\lim_{t\ra \un}\frac d{dt}f(u)=-\tfrac1n.
\end{equation}

\cvm
The same result is of course valid for $\frac d{dt}\tilde f(\tilde u)$, where one should
note the ambiguous usage of the tilde.

\cvm
Combining \re{3.4}, \re{3.6} and de L'Hospital's rule we infer that $\om=\tilde \om$.

\cvm
To prove
\begin{equation}
\lim \int_{M(t)}G_{\al\bet}\nu^\al\nu^\bet e^{\om f} e^\psi= \lim
\int_{M(t)}G_{\al\bet}\nu^\al\nu^\bet e^{\om \tilde f} e^{\tilde\psi}
\end{equation}
it suffices to show that $m=\tilde m$, or equivalently,
\begin{equation}
\lim \frac{e^{\tilde f}}{e^f}=1.
\end{equation}

Since $\om=\tilde \om$, the relation \re{3.4} implies
\begin{equation}
\lim \frac{e^{\tilde f}}{e^f}=c=\const.
\end{equation}

Let $(\x^i)$ be local coordinates for $M(t)$ and  $x=x(\x)$ be a local embedding,
then
\begin{equation}
g_{ij}=\spd{x_i}{x_j}
\end{equation}
is the induced metric. Let
\begin{equation}\lae{3.11}
\tilde g_{ij}=e^{-2(f+\psi)} g_{ij}
\end{equation}
be the conformal metric, then
\begin{equation}
\int_\so \sqrt{\det(\tilde g_{ij})}=\int_\so v \sqrt{\det(\s_{ij}(u,x)},
\end{equation}
where
\begin{equation}
v^2=1-\abs{Du}^2=1-\s^{ij}u_iu_j
\end{equation}
and hence
\begin{equation}
\lim \int_\so v \sqrt{\det(\s_{ij}(u,x))}=\abs{S^n}
\end{equation}
due to our normalization assumption.

\cvm
On the other hand, if we express the right-hand side of \re{3.11} with respect to the
second coordinate system, then we obtain
\begin{equation}
\abs{S^n}=\lim \int_\so e^{n(\tilde f+\tilde \psi)}e^{-n(f+\psi)}v
\sqrt{\det(\tilde\s_{ij}(\tilde u,x))}=c^n \abs{S^n},
\end{equation}
hence $c=1$.

\section{A variant of the Penrose inequality}

Let $M_\tau=\{x^0=\tau\}$ be coordinate slices and suppose that the integrals
\begin{equation}\lae{4.1}
\int_{M_\tau}G_{\al\bet}\nu^\al\nu^\bet e^{\om f}e^\psi
\end{equation}
would increase monotonically with respect to $\tau$, then
\begin{equation}
\lim_{\tau\ra 0}\int_{M_\tau}G_{\al\bet}\nu^\al\nu^\bet e^{\om f}e^\psi\ge
\int_{M_{\bar\tau}}G_{\al\bet}\nu^\al\nu^\bet e^{\om f}e^\psi.
\end{equation}
If $M_{\bar\tau}$ would be totally geodesic, then
\begin{equation}
\int_{M_{\bar\tau}}G_{\al\bet}\nu^\al\nu^\bet e^{\om f}e^\psi=\tfrac12
\int_{M_{\bar\tau}}R e^{\om f}e^\psi
\end{equation}

\cvm
To prove the monotonicity of the integrals in \re{4.1}, let us look at the relation
\re{2.11}. The monotonicity of the integrals is equivalent to  the non-negativity of
the right-hand side of \re{2.11}. This will be the case for highly symmetrical
spacetimes as we shall see in the next section. For general ARW spacetimes however,
non-negativity of the right-hand side of \re{2.11} could only be derived under the
assumptions $\om=0$, $\psi=0$, and, furthermore, that the slices $M_\tau$ are
convex, i.e., $\bar h_{ij}\ge 0$, and the spatial part of the Einstein tensor positive
semi-definite, i.e., $G^{ij}\ge 0$.

Notice that $f'$ is always negative and that, under physical assumptions, $G^{00}\ge
0$ and also that the signs of the spatial part $(G^{ij})$ and of $\om$ are the same,
as can be derived from the Einstein equations
\begin{equation}
G_{\al\bet}=\ka T_{\al\bet},
\end{equation}
if the stress energy tensor is supposed to be asymptotically that of a perfect fluid
with an equation of state
\begin{equation}
p=\tfrac\om n \rho.
\end{equation}

\section{An example}

Let $\hat N$ be the $\tup{S-AdS}_{(n+2)}$ spacetime with metric
\begin{equation}
d\hat s^2=-h dt^2+ h^{-1}dr^2+r^2 \s_{ij}dx^idx^j,
\end{equation}
where
\begin{equation}
h=1-\frac2{n(n+1)}\Lam r^2-m r^{-(n-1)}
\end{equation}
with constants $\Lam\le 0$ and $m>0$; $(\s_{ij})$ is the metric of $S^n$.

In $r=0$ there is a black hole singularity, the event horizon is in $r=r_0$, such that
$h(r_0)=0$, and the black hole region is given by $\{h<0\}=\{0<r<r_0\}$.

\cvm
In the black hole region $t$ is a spatial coordinate and $r$ the time coordinate. Set
\begin{equation}
\tilde h=- h
\end{equation}
and consider in the black hole region the brane
\begin{equation}
N=\{t=\const, 0<r<r_0\}.
\end{equation}

The induced metric $(\bar g_{\al\bet})$ is
\begin{equation}
\begin{aligned}
d\bar s^2&=-\tilde h^{-1} dr^2 +r^2 \s_{ij}dx^idx^j\\
&=r^2(-r^{-2}\tilde h^{-1} dr^2+\s_{ij}dx^idx^j).
\end{aligned}
\end{equation}

Define
\begin{equation}
f=\log r
\end{equation}
and $x^0$ by
\begin{equation}
dx^0=-r^{-1}\tilde h^{-\frac12} dr,
\end{equation}
i.e.,
\begin{equation}
x^0=-\int_0^rs^{-1}\tilde h^{-\frac12} ds.
\end{equation}

$x^0$ is then a future directed time function and the (induced) singularity lies in
$x^0=0$.

In these coordinates the metric has the form
\begin{equation}
d\bar s^2=e^{2f}(-(dx^0)^2+\s_{ij}dx^idx^j).
\end{equation}

Let a prime denote differentiation with respect to $x^0$ and a dot with respect to
$r$, then
\begin{equation}
f'=r^{-1}\frac{dr}{dx^0}=-\tilde h^{\frac12}
\end{equation}
and
\begin{equation}
f''=\tfrac12 r \dot{\tilde h}=-\tfrac12 m(n-1)r^{-(n-1)}+\frac2{n(n+1)}\Lam.
\end{equation}

\cvm
Set $\om =1$ so that
\begin{equation}
\tilde\ga=\tfrac12(n+\om-2)=\tfrac12(n-1),
\end{equation}
then
\begin{equation}
f''+\tilde\ga \abs{f'}^2=\tfrac1n\Lam r^2-\tfrac12 (n-1)
\end{equation}
and
\begin{equation}
\abs{f'}^2e^{2\tilde\ga f}=\abs{f'}^2 r^{(n-1)}=\tilde h
r^{(n-1)}=m+\frac2{n(n+1)}\Lam r^{(n+1)}-r^{(n-1)}.
\end{equation}

\cvm
Let $M_r$ be the coordinate slices $\{x^0=\const\}$, then the past directed normal
is $\nu=(\nu^\al)=e^{-f}(-1,\ldots,0)$,
\begin{equation}
G_{\al\bet}\nu^\al\nu^\bet  e^{2 f}=\tfrac12 n(n-1)\abs{f'}^2+\tfrac12 n(n-1)=\tfrac12
n(n-1) (\tilde h+1),
\end{equation}
and
\begin{equation}
\int_{M_r}G_{\al\bet}\nu^\al\nu^\bet e^f=\tfrac12 n(n-1) \int_{S^n}(\tilde
h+1)r^{(n-1)},
\end{equation}
hence
\begin{equation}
\lim_{r\ra 0}\int_{M_r}G_{\al\bet}\nu^\al\nu^\bet e^f =\tfrac12 n(n-1) \abs{S^n}
m.
\end{equation}

The convergence is monotone increasing in $x^0$.

\cvm
Note that the value $\om=1$ corresponds to a radiation dominated universe if the
stress energy tensor is asymptotically equal to that of a perfect fluid and the
equation of state is
\begin{equation}
p=\tfrac\om n\rho.
\end{equation}
\nocite{cg:imcf}
\bibliographystyle{hamsplain}
\bibliography{mrabbrev,publications}
\providecommand{\bysame}{\leavevmode\hbox to3em{\hrulefill}\thinspace}
\providecommand{\href}[2]{#2}


\end{document}